%version of 12.3.1996 with changes made by MF
\magnification=\magstep1%metka
\input amstex
\def\temp{1.35}%
\let\tempp=\relax
\expandafter\ifx\csname psboxversion\endcsname\relax
  \message{PSBOX(\temp)}%
\else
    \ifdim\temp cm>\psboxversion cm
      \message{PSBOX(\temp)}%
    \else
      \message{PSBOX(\psboxversion) is already loaded: I won't load
        PSBOX(\temp)!}%
      \let\temp=\psboxversion
      \let\tempp= 
    \fi
\fi
\tempp
\message{by Jean Orloff: loading ...}
\let\psboxversion=\temp
\catcode`\@=11
% Every macro likes a little privacy...
%
%Trying to tame the variety of \special commands for Postscript: the
%  universal internal command \PSspeci@l##1##2 takes ##1 to be the
%  filename and ##2 to be the integer scale factor*1000 (as for usual
%   TeX \scale commands)
%
\def\psfortextures{%     For TeXtures on the Macintosh
%-----------------
\def\PSspeci@l##1##2{%
\special{illustration ##1\space scaled ##2}%
}}%
\def\psfordvitops{%      For the DVItoPS converter on IBM mainframes
%----------------
\def\PSspeci@l##1##2{%
\special{dvitops: import ##1\space \the\drawingwd \the\drawinght}%
}}%
\def\psfordvips{%      For DVIPS converter on VAX, UNIX and PC's
%--------------
\def\PSspeci@l##1##2{%
%    \special{/@scaleunit 1000 def}% never read dox without trying!
\d@my=0.1bp \d@mx=\drawingwd \divide\d@mx by\d@my% BUG! for large \drawingwd
\includegraphics{##1\space}}}%
\def\psforoztex{%        For the OzTeX shareware on the Macintosh
%--------------
\def\PSspeci@l##1##2{%
\special{##1 \space
      ##2 1000 div dup scale
      \number-\psllx\space\space \number-\pslly\space\space translate
}}}%
\def\psfordvitps{%       From the UNIX TeXPS package, vers.>3.12
%---------------
\def\dvitpsLiter@ldim##1{\dimen0=##1\relax
\special{dvitps: Literal "\number\dimen0\space"}}%
\def\PSspeci@l##1##2{%
% psfig.psr contains the def of "startTexFig": if you can locate it
% and put it somewhere in the TEXINPUTS path, this works
\at(0bp;\drawinght){%
\special{dvitps: Include0 "psfig.psr"}% contains def of "startTexFig"
\dvitpsLiter@ldim{\drawingwd}%
\dvitpsLiter@ldim{\drawinght}%
\dvitpsLiter@ldim{\psllx bp}%
\dvitpsLiter@ldim{\pslly bp}%
\dvitpsLiter@ldim{\psurx bp}%
\dvitpsLiter@ldim{\psury bp}%
\special{dvitps: Literal "startTexFig"}%
\special{dvitps: Include1 "##1"}%
\special{dvitps: Literal "endTexFig"}%
}}}%
\def\psfordvialw{%   Try for dvialw, a UNIX public domain
%---------------
\def\PSspeci@l##1##2{
\special{language "PostScript",
position = "bottom left",
literal "  \psllx\space \pslly\space translate
  ##2 1000 div dup scale
  -\psllx\space -\pslly\space translate",
include "##1"}
}}%
\def\psforptips{%   For MS-DOS; LUOMA@brandeis.bitnet
%---------------
\def\PSspeci@l##1##2{{
\d@mx=\psurx bp
\advance \d@mx by -\psllx bp
\divide \d@mx by 1000\multiply\d@mx by \xscale
\incm{\d@mx}
\let\tmpx\dimincm
\d@my=\psury bp
\advance \d@my by -\pslly bp
\divide \d@my by 1000\multiply\d@my by \xscale
\incm{\d@my}
\let\tmpy\dimincm
\d@mx=-\psllx bp
\divide \d@mx by 1000\multiply\d@mx by \xscale
\d@my=-\pslly bp
\divide \d@my by 1000\multiply\d@my by \xscale
\at(\d@mx;\d@my){\special{ps:##1 x=\tmpx cm, y=\tmpy cm}}
}}}%
\def\psonlyboxes{%     Draft-like behaviour if none of the others works
%---------------
\def\PSspeci@l##1##2{%
\at(0cm;0cm){\boxit{\vbox to\drawinght
  {\vss\hbox to\drawingwd{\at(0cm;0cm){\hbox{({\tt##1})}}\hss}}}}
}}%
\def\psloc@lerr#1{%
\let\savedPSspeci@l=\PSspeci@l%
\def\PSspeci@l##1##2{%
\at(0cm;0cm){\boxit{\vbox to\drawinght
  {\vss\hbox to\drawingwd{\at(0cm;0cm){\hbox{({\tt##1}) #1}}\hss}}}}
\let\PSspeci@l=\savedPSspeci@l% restore normal output for other figs!
}}%
%\def\psfor...  add your own!
%
% Some common defs
%
\newread\pst@mpin
\newdimen\drawinght\newdimen\drawingwd
\newdimen\psxoffset\newdimen\psyoffset
\newbox\drawingBox
\newcount\xscale \newcount\yscale \newdimen\pscm\pscm=1cm
\newdimen\d@mx \newdimen\d@my
\newdimen\pswdincr \newdimen\pshtincr
\let\ps@nnotation=\relax
{\catcode`\|=0 |catcode`|\=12 |catcode`|%=12 |catcode`~=12
|catcode`#=12 |catcode`*=14
|xdef|backslashother{\}*
|xdef|percentother{%}*
|xdef|tildeother{~}*
|xdef|sharpother{#}*
}%
% useful to display special chars in \tt; fails for \,#,%
\def\R@moveMeaningHeader#1:->{}%
\def\uncatcode#1{%
\edef#1{\expandafter\R@moveMeaningHeader\meaning#1}}%
\def\execute#1{#1}% NOT stupid: cs in #1 are then identified BEFORE execution
\def\psm@keother#1{\catcode`#112\relax}% borrowed from latex
\def\executeinspecs#1{%
\execute{\begingroup\let\do\psm@keother\dospecials\catcode`\^^M=9#1\endgroup}}%
\def\@mpty{}%
% \if\matchin#1#2<=> \iftrue if #1 contains #2, <=>\iffalse otherwise:
% \if\matchexpin: idem, but #1 & #2 are first fully expanded (no \if
% inside!)
% \tmpa & \tmpb contain what's before and after the occurence of #2
\def\matchexpin#1#2{
  \fi%
%\message{(#1>#2)}
  \edef\tmpb{{#2}}%
  \expandafter\makem@tchtmp\tmpb%
  \edef\tmpa{#1}\edef\tmpb{#2}%
  \expandafter\expandafter\expandafter\m@tchtmp\expandafter\tmpa\tmpb\endm@tch%
  \if\match%
}%
\def\matchin#1#2{%
  \fi%
  \makem@tchtmp{#2}%
  \m@tchtmp#1#2\endm@tch%
  \if\match%
}%
\def\makem@tchtmp#1{\def\m@tchtmp##1#1##2\endm@tch{%
  \def\tmpa{##1}\def\tmpb{##2}\let\m@tchtmp=\relax%
  \ifx\tmpb\@mpty\def\match{YN}%
  \else\def\match{YY}\fi%
}}%
% converts any dimen in cm, with 1E-4 cm precision
\def\incm#1{{\psxoffset=1cm\d@my=#1
 \d@mx=\d@my
  \divide\d@mx by \psxoffset
  \xdef\dimincm{\number\d@mx.}
  \advance\d@my by -\number\d@mx cm
  \multiply\d@my by 100
 \d@mx=\d@my
  \divide\d@mx by \psxoffset
  \edef\dimincm{\dimincm\number\d@mx}
  \advance\d@my by -\number\d@mx cm
  \multiply\d@my by 100
 \d@mx=\d@my
  \divide\d@mx by \psxoffset
  \xdef\dimincm{\dimincm\number\d@mx}
}}%
%
%  \ReadPSize{PSfilename} reads the dimensions of a PostScript drawing
%      and stores it in \drawinght(wd)
\newif\ifNotB@undingBox
\newhelp\PShelp{Proceed: you'll have a 5cm square blank box instead of
your graphics.}%
\def\s@tsize#1 #2 #3 #4\@ndsize{
  \def\psllx{#1}\def\pslly{#2}%
  \def\psurx{#3}\def\psury{#4}%  needed by a crazyness of dvips!
  \ifx\psurx\@mpty\NotB@undingBoxtrue% this is not a valid one!
  \else
    \drawinght=#4bp\advance\drawinght by-#2bp
    \drawingwd=#3bp\advance\drawingwd by-#1bp
%  !Units related by crazy factors as bp/pt=72.27/72 should be BANNED!
  \fi
  }%
\def\sc@nBBline#1:#2\@ndBBline{\edef\p@rameter{#1}\edef\v@lue{#2}}%
\def\g@bblefirstblank#1#2:{\ifx#1 \else#1\fi#2}%
{\catcode`\%=12
\xdef\B@undingBox{%%BoundingBox}}%
%% is not a true comment in PostScript, even if % is!
\def\ReadPSize#1{
 \readfilename#1\relax
 \let\PSfilename=\lastreadfilename
 \openin\pst@mpin=#1\relax
 \ifeof\pst@mpin \errhelp=\PShelp
   \errmessage{I haven't found your postscript file (\PSfilename)}%
   \psloc@lerr{was not found}%
   \s@tsize 0 0 142 142\@ndsize
   \closein\pst@mpin
 \else
% each entry in \GlobalInputList should be unique
   \if\matchexpin{\GlobalInputList}{, \lastreadfilename}%
   \else\xdef\GlobalInputList{\GlobalInputList, \lastreadfilename}%
     \immediate\write\psbj@inaux{\lastreadfilename,}%
   \fi%
   \loop
     \executeinspecs{\catcode`\ =10\global\read\pst@mpin to\n@xtline}%
     \ifeof\pst@mpin
       \errhelp=\PShelp
       \errmessage{(\PSfilename) is not an Encapsulated PostScript File:
           I could not find any \B@undingBox: line.}%
       \edef\v@lue{0 0 142 142:}%
       \psloc@lerr{is not an EPSFile}%
       \NotB@undingBoxfalse
     \else
       \expandafter\sc@nBBline\n@xtline:\@ndBBline
       \ifx\p@rameter\B@undingBox\NotB@undingBoxfalse
         \edef\t@mp{%
           \expandafter\g@bblefirstblank\v@lue\space\space\space}%
         \expandafter\s@tsize\t@mp\@ndsize
       \else\NotB@undingBoxtrue
       \fi
     \fi
   \ifNotB@undingBox\repeat
   \closein\pst@mpin
 \fi
\message{#1}%
}%
%
% \psboxto(xdim;ydim){psfilename}: you specify the dimensions and
%    TeX uniformly scales to fit the largest one. If xdim=0pt, the
%    scale is fully determined by ydim and vice versa.
%    Notice: psboxes are a real vboxes; couldn't take hbox otherwise all
%    indentation and all cr's would be interpreted as spaces (hugh!).
%
\def\psboxto(#1;#2)#3{\vbox{%
   \ReadPSize{#3}%
   \advance\pswdincr by \drawingwd
   \advance\pshtincr by \drawinght
   \divide\pswdincr by 1000
   \divide\pshtincr by 1000
   \d@mx=#1
   \ifdim\d@mx=0pt\xscale=1000
         \else \xscale=\d@mx \divide \xscale by \pswdincr\fi
   \d@my=#2
   \ifdim\d@my=0pt\yscale=1000
         \else \yscale=\d@my \divide \yscale by \pshtincr\fi
   \ifnum\yscale=1000
         \else\ifnum\xscale=1000\xscale=\yscale
                    \else\ifnum\yscale<\xscale\xscale=\yscale\fi
              \fi
   \fi
   \divide\drawingwd by1000 \multiply\drawingwd by\xscale
   \divide\drawinght by1000 \multiply\drawinght by\xscale
   \divide\psxoffset by1000 \multiply\psxoffset by\xscale
   \divide\psyoffset by1000 \multiply\psyoffset by\xscale
   \global\divide\pscm by 1000
   \global\multiply\pscm by\xscale
   \multiply\pswdincr by\xscale \multiply\pshtincr by\xscale
   \ifdim\d@mx=0pt\d@mx=\pswdincr\fi
   \ifdim\d@my=0pt\d@my=\pshtincr\fi
   \message{scaled \the\xscale}%
 \hbox to\d@mx{\hss\vbox to\d@my{\vss
   \global\setbox\drawingBox=\hbox to 0pt{\kern\psxoffset\vbox to 0pt{%
      \kern-\psyoffset
      \PSspeci@l{\PSfilename}{\the\xscale}%
      \vss}\hss\ps@nnotation}%
   \global\wd\drawingBox=\the\pswdincr
   \global\ht\drawingBox=\the\pshtincr
   \global\drawingwd=\pswdincr
   \global\drawinght=\pshtincr
   \baselineskip=0pt
   \copy\drawingBox
 \vss}\hss}%
  \global\psxoffset=0pt
  \global\psyoffset=0pt
  \global\pswdincr=0pt
  \global\pshtincr=0pt % These are local to one figure
  \global\pscm=1cm %should not be necessary
}}%
%
% \psboxscaled{scalefactor*1000}{PSfilename} allows to bypass the
%   rounding errors of TeX integer divisions for situations where the
%   TeX box should fit the original BoundingBox with a precision
%   better
%   than 1/1000.
%
\def\psboxscaled#1#2{\vbox{%
  \ReadPSize{#2}%
  \xscale=#1
  \message{scaled \the\xscale}%
  \divide\pswdincr by 1000 \multiply\pswdincr by \xscale
  \divide\pshtincr by 1000 \multiply\pshtincr by \xscale
  \divide\psxoffset by1000 \multiply\psxoffset by\xscale
  \divide\psyoffset by1000 \multiply\psyoffset by\xscale
  \divide\drawingwd by1000 \multiply\drawingwd by\xscale
  \divide\drawinght by1000 \multiply\drawinght by\xscale
  \global\divide\pscm by 1000
  \global\multiply\pscm by\xscale
  \global\setbox\drawingBox=\hbox to 0pt{\kern\psxoffset\vbox to 0pt{%
     \kern-\psyoffset
     \PSspeci@l{\PSfilename}{\the\xscale}%
     \vss}\hss\ps@nnotation}%
  \advance\pswdincr by \drawingwd
  \advance\pshtincr by \drawinght
  \global\wd\drawingBox=\the\pswdincr
  \global\ht\drawingBox=\the\pshtincr
  \global\drawingwd=\pswdincr
  \global\drawinght=\pshtincr
  \baselineskip=0pt
  \copy\drawingBox
  \global\psxoffset=0pt
  \global\psyoffset=0pt
  \global\pswdincr=0pt
  \global\pshtincr=0pt % These are local to one figure
  \global\pscm=1cm
}}%
%
%  \psbox{PSfilename} makes a TeX box having the minimal size to
%      enclose the picture
\def\psbox#1{\psboxscaled{1000}{#1}}%
%------------------------------------------------------
%  \joinfiles file1, file2, ...n \into joinedfilename .
%     makes one file out of many
%  \splitfile joinedfilename
%     the opposite
\newif\ifn@teof\n@teoftrue
\newif\ifc@ntrolline
\newif\ifmatch
\newread\j@insplitin
\newwrite\j@insplitout
\newwrite\psbj@inaux
\immediate\openout\psbj@inaux=psbjoin.aux
\immediate\write\psbj@inaux{\string\joinfiles}%
\immediate\write\psbj@inaux{\jobname,}%
%
% INPUT REDEFINITION
%
% works if #1 is a single character
\def\toother#1{\ifcat\relax#1\else\expandafter%
  \toother@ux\meaning#1\endtoother@ux\fi}%
\def\toother@ux#1 #2#3\endtoother@ux{\def\tmp{#3}%
  \ifx\tmp\@mpty\def\tmp{#2}\let\next=\relax%
  \else\def\next{\toother@ux#2#3\endtoother@ux}\fi%
\next}%
%
% \readfilename defs:
%
\let\readfilenamehook=\relax
\def\re@d{\expandafter\re@daux}% spares typing 10 \expandafter's...
\def\re@daux{\futurelet\nextchar\stopre@dtest}%
\def\re@dnext{\xdef\lastreadfilename{\lastreadfilename\nextchar}%
  \afterassignment\re@d\let\nextchar}%
\def\stopre@d{\egroup\readfilenamehook}%
\def\stopre@dtest{%
  \ifcat\nextchar\relax\let\nextread\stopre@d
  \else
    \ifcat\nextchar\space\def\nextread{%
      \afterassignment\stopre@d\chardef\nextchar=`}%
    \else\let\nextread=\re@dnext
      \toother\nextchar
      \edef\nextchar{\tmp}%
    \fi
  \fi\nextread}%
\def\readfilename{\bgroup%
  \let\\=\backslashother \let\%=\percentother \let\~=\tildeother
  \let\#=\sharpother \xdef\lastreadfilename{}%
  \re@d}%
%
% redefines \input using \readfilename
%
\xdef\GlobalInputList{\jobname}%
\def\psnewinput{%
  \def\readfilenamehook{% each entry in \GlobalInputList should be unique
    \if\matchexpin{\GlobalInputList}{, \lastreadfilename}%
    \else\xdef\GlobalInputList{\GlobalInputList, \lastreadfilename}%
      \immediate\write\psbj@inaux{\lastreadfilename,}%
    \fi%
    \let\readfilenamehook=\relax%
    \ps@ldinput\lastreadfilename\relax%
  }\readfilename%
}%
\expandafter\ifx\csname @@input\endcsname\relax    % then Plain
  \immediate\let\ps@ldinput=\input\def\input{\psnewinput}%
\else
  \immediate\let\ps@ldinput=\@@input
  \def\@@input{\psnewinput}%
\fi%
\def\nowarnopenout{%
 \def\warnopenout##1##2{%
   \readfilename##2\relax
   \message{\lastreadfilename}%
   \immediate\openout##1=\lastreadfilename\relax}}%
\def\warnopenout#1#2{%
 \readfilename#2\relax
 \def\t@mp{TrashMe,psbjoin.aux,psbjoint.tex,}\uncatcode\t@mp
 \if\matchexpin{\t@mp}{\lastreadfilename,}%
 \else
   \immediate\openin\pst@mpin=\lastreadfilename\relax
   \ifeof\pst@mpin
     \else
     \edef\tmp{{If the content of this file is precious to you, this
is your last chance to abort (ie press x or e) and rename it before
retexing (\jobname). If you're sure there's no file
(\lastreadfilename) in the directory of (\jobname), then go on: I'm
simply worried because you have another (\lastreadfilename) in some
directory I'm looking in for inputs...}}%
     \errhelp=\tmp
     \errmessage{I may be about to replace your file named \lastreadfilename}%
   \fi
   \immediate\closein\pst@mpin
 \fi
 \message{\lastreadfilename}%
 \immediate\openout#1=\lastreadfilename\relax}%
% % will have an unusual catcode below; use * instead
{\catcode`\%=12\catcode`\*=14
\gdef\splitfile#1{*
 \readfilename#1\relax
 \immediate\openin\j@insplitin=\lastreadfilename\relax
 \ifeof\j@insplitin
   \message{! I couldn't find and split \lastreadfilename!}*
 \else
   \immediate\openout\j@insplitout=TrashMe
   \message{< Splitting \lastreadfilename\space into}*
   \loop
     \ifeof\j@insplitin
       \immediate\closein\j@insplitin\n@teoffalse
     \else
       \n@teoftrue
       \executeinspecs{\global\read\j@insplitin to\spl@tinline\expandafter
         \ch@ckbeginnewfile\spl@tinline%Beginning-Of-File-Named:%\endcheck}*
       \ifc@ntrolline
       \else
         \toks0=\expandafter{\spl@tinline}*
         \immediate\write\j@insplitout{\the\toks0}*
       \fi
     \fi
   \ifn@teof\repeat
   \immediate\closeout\j@insplitout
 \fi\message{>}*
}*
\gdef\ch@ckbeginnewfile#1%Beginning-Of-File-Named:#2%#3\endcheck{*
 \def\t@mp{#1}*
 \ifx\@mpty\t@mp
   \def\t@mp{#3}*
   \ifx\@mpty\t@mp
     \global\c@ntrollinefalse
   \else
     \immediate\closeout\j@insplitout
     \warnopenout\j@insplitout{#2}*
     \global\c@ntrollinetrue
   \fi
 \else
   \global\c@ntrollinefalse
 \fi}*
\gdef\joinfiles#1\into#2{*
 \message{< Joining following files into}*
 \warnopenout\j@insplitout{#2}*
 \message{:}*
 {*
 \edef\w@##1{\immediate\write\j@insplitout{##1}}*
\w@{% This collection of files was produced with CERN psbox package}*
\w@{% To decompose and tex it:}*
\w@{%-save this with a filename CONTAINING ONLY LETTERS and a .TEX}*
\w@{% extension (say, JOINTFIL.TEX), in some empty directory;}*
\w@{%-make sure you can \string\input\space psbox.tex (version>=1.3);}*
\w@{%  (else ftp cs.nyu.edu(=128.122.140.24):pub/TeX/psbox/, then get}*
\w@{%  and tex the file psboxall.tex; more info in psbREAD.ME)}*
\w@{%-tex JOINTFIL.TEX using Plain, or LaTeX, or whatever is needed by}*
\w@{%  the first file in the joining (after splitting JOINTFIL.TEX into}*
\w@{%  it's constituents, TeX will try to process it as it stands).}*
\w@{\string\input\space psbox.tex}*
\w@{\string\splitfile{\string\jobname}}*
\w@{\string\let\string\autojoin=\string\relax}*
}*
 \expandafter\tre@tfilelist#1, \endtre@t
 \immediate\closeout\j@insplitout
 \message{>}*
}*
\gdef\tre@tfilelist#1, #2\endtre@t{*
 \readfilename#1\relax
 \ifx\@mpty\lastreadfilename
 \else
   \immediate\openin\j@insplitin=\lastreadfilename\relax
   \ifeof\j@insplitin
     \errmessage{I couldn't find file \lastreadfilename}*
   \else
     \message{\lastreadfilename}*
     \immediate\write\j@insplitout{%Beginning-Of-File-Named:\lastreadfilename}*
     \executeinspecs{\global\read\j@insplitin to\oldj@ininline}*
     \loop
       \ifeof\j@insplitin\immediate\closein\j@insplitin\n@teoffalse
       \else\n@teoftrue
         \executeinspecs{\global\read\j@insplitin to\j@ininline}*
         \toks0=\expandafter{\oldj@ininline}*
         \let\oldj@ininline=\j@ininline
         \immediate\write\j@insplitout{\the\toks0}*
       \fi
     \ifn@teof
     \repeat
   \immediate\closein\j@insplitin
   \fi
   \tre@tfilelist#2, \endtre@t
 \fi}*
}%
% To be put at the end of a file, for making a tar-like file containing
%   everything it used.
\def\autojoin{%
 \immediate\write\psbj@inaux{\string\into{psbjoint.tex}}%
 \immediate\closeout\psbj@inaux
 \expandafter\joinfiles\GlobalInputList\into{psbjoint.tex}%
}%
%----------------------------------------------------------------
%  Annotations & Captions etc...
%
%
% \centinsert{anybox} is just a centered \midinsert, but is included as
%    people barely use the original inserts from TeX.
%
\def\centinsert#1{\midinsert\line{\hss#1\hss}\endinsert}%
\def\psannotate#1#2{\vbox{%
  \def\ps@nnotation{#2\global\let\ps@nnotation=\relax}#1}}%
\def\pscaption#1#2{\vbox{%
   \setbox\drawingBox=#1
   \copy\drawingBox
   \vskip\baselineskip
   \vbox{\hsize=\wd\drawingBox\setbox0=\hbox{#2}%
     \ifdim\wd0>\hsize
       \noindent\unhbox0\tolerance=5000
    \else\centerline{\box0}%
    \fi
}}}%
% for compatibility with older versions, but \psfig is a bad name!
%\def\psfig#1#2#3{\pscaption{\psannotate{#1}{#2}}{#3}}
%\def\psfigurebox#1#2#3{\pscaption{\psannotate{\psbox{#1}}{#2}}{#3}}
%
% \at(#1;#2)#3 puts #3 at #1-higher and #2-right of the current
%    position without moving it (to be used in annotations).
\def\at(#1;#2)#3{\setbox0=\hbox{#3}\ht0=0pt\dp0=0pt
  \rlap{\kern#1\vbox to0pt{\kern-#2\box0\vss}}}%
%
% \gridfill(ht;wd) makes a 1cm*1cm grid of ht by wd whose lower-left
%   corner is the current point
\newdimen\gridht \newdimen\gridwd
\def\gridfill(#1;#2){%
  \setbox0=\hbox to 1\pscm
  {\vrule height1\pscm width.4pt\leaders\hrule\hfill}%
  \gridht=#1
  \divide\gridht by \ht0
  \multiply\gridht by \ht0
  \gridwd=#2
  \divide\gridwd by \wd0
  \multiply\gridwd by \wd0
  \advance \gridwd by \wd0
  \vbox to \gridht{\leaders\hbox to\gridwd{\leaders\box0\hfill}\vfill}}%
%
% Useful to measure where to put annotations
\def\fillinggrid{\at(0cm;0cm){\vbox{%
  \gridfill(\drawinght;\drawingwd)}}}%
%
% \textleftof\anybox: Sample text\endtext
%   inserts "Sample text" on the left of \anybox ie \vbox, \psbox.
%   \textrightof is the symmetric (not documented, too uggly)
% Welcome any suggestion about clean wraparound macros from
%   TeXhackers reading this
%
\def\textleftof#1:{%
  \setbox1=#1
  \setbox0=\vbox\bgroup
    \advance\hsize by -\wd1 \advance\hsize by -2em}%
\def\textrightof#1:{%
  \setbox0=#1
  \setbox1=\vbox\bgroup
    \advance\hsize by -\wd0 \advance\hsize by -2em}%
\def\endtext{%
  \egroup
  \hbox to \hsize{\valign{\vfil##\vfil\cr%
\box0\cr%
\noalign{\hss}\box1\cr}}}%
%
% \frameit{\thick}{\skip}{\anybox}
%    draws with thickness \thick a box around \anybox, leaving \skip of
%    blank around it. eg \frameit{0.5pt}{1pt}{\hbox{hello}}
% \boxit{\anybox} is a shortcut.
\def\frameit#1#2#3{\hbox{\vrule width#1\vbox{%
  \hrule height#1\vskip#2\hbox{\hskip#2\vbox{#3}\hskip#2}%
        \vskip#2\hrule height#1}\vrule width#1}}%
\def\boxit#1{\frameit{0.4pt}{0pt}{#1}}%
\catcode`\@=12 % cs containing @ are unreachable
%
% CUSTOMIZE YOUR DEFAULT DRIVER:
%    Uncomment the line corresponding to your TeX system:
% \psfortextures%     For TeXtures on the Macintosh
%\psforoztex   %     For OzTeX shareware on the Macintosh
%\psfordvitops %     For the DVItoPS converter for TeX on IBM mainframes
\psfordvips   %     For DVIPS converter on VAX and UNIX
%\psfordvitps  %     For dvitps from TeXPS package under UNIX
%\psfordvialw  %     For dvialw, UNIX public domain
%\psonlyboxes  %     Blank Boxes (when all else fails).

\documentstyle{amsppt}
\TagsOnRight
\hsize=5in                                                  
\vsize=7.8in
\strut
\vskip5truemm
\font\small=cmr8
\font\itsmall=cmti8
%\font\ref=cmr9
%\font\refit=cmti9
%\font\refbf=cmbx9

\def\smallarea#1{\par\begingroup\baselineskip=10pt#1\endgroup
\par}

\def\refer#1#2{\par\begingroup\baselineskip=11pt\noindent
\leftskip=8.25mm\rightskip=0mm

\strut\llap{#1\kern 1em}{#2\hfill}\par\endgroup}
	      
\nopagenumbers
%\subheading{Revised version}
%\vskip 2cm

\topmatter
\title Morse theory of harmonic forms\endtitle
%\rightheadtext{}
\leftheadtext{M.Farber, G.Katz and J.Levine}
\author  Michael Farber, Gabriel Katz and Jerome Levine \endauthor
\address
School of Mathematical Sciences, Tel-Aviv University, Ramat-Aviv 
69978,
 Israel
\endaddress
\address Department of Mathematics, Clark University, Worcester, 
MA \endaddress
\address Department of Mathematics, Brandeis University, Waltham, 
MA 
\endaddress
\email farber\@math.tau.ac.il\quad  
levine\@math.binah.cc.brandeis.edu
\quad  gkatz\@bbn.com 
\endemail
\thanks{The research was  supported by US - Israel
Binational Science Foundation Grants 9400299 and 9400073, and by 
NSF 
Grant
 93-03489.}
\endthanks
\abstract{ We consider the problem of whether it is possible 
to improve the
Novikov inequalities for closed 1-forms, or any other inequalities of 
a 
similar nature, if we assume, additionally, that the given 1-form is 
harmonic with
respect to some Riemannian metric. We show that, under suitable 
assumptions, it is impossible. We use a theorem of E.Calabi \cite{C}, 
characterizing 1-forms which are harmonic with respect to some 
metric, in an
essential way. We also study some interesting examples illustrating 
our results.}
\endabstract 
\endtopmatter

\define\C{{\Bbb C}}
\define\R{{\Bbb R}}  
\define\RR{{\Cal R}}       
\define\Z{{\Bbb Z}}  
\define\Q{{\Bbb Q}} 
\define\V{{\tilde V}}
\define\OO{{\Cal O}}  
\define\T{{\Cal T}}  
\define\fs{{\Cal S}}
\define\N{{\Bbb N}}
\define\Res{\operatorname{Res}} 
\define\Hom{\operatorname{Hom}}  
\define\Herm{\operatorname{Herm}}   
\define\ch{\operatorname{ch}} 
\define\Tr{\operatorname{Tr}}   
\define\Tor{\operatorname{Tor}} 
\define\Det{\operatorname{Det}}   
\define\Ext{\operatorname{Ext}} 
\define\Aut{\operatorname{Aut}}
\define\rk{\operatorname{rk}}          
\define\index{\operatorname{index}}  
\define\Ker{\operatorname{Ker}} 
\define\End{\operatorname{End}}   
\define\Har{\operatorname{Har}}
\define\GL{\operatorname{GL}} 
\define\Diff{\operatorname{Diff}}        
\define\M{{\Cal M}} 
\define\A{{\Cal A}}  
\redefine\B{{\Cal B}}  
\define\p{{\frak p}} 
\define\Cl{{\frak Cl}}    
\define\F{{\Cal F}}      
\define\sign{\operatorname{sign}} 
\define\tr{\operatorname{tr}} 
\define\im{\operatorname{im}} 
\define\id{\operatorname{id}}    
\define\cl{\operatorname{cl}}
\define\Int{\operatorname{Int}}  
%\define\mod{\operatorname{mod}}
\redefine\H{\Cal H}
\redefine\D{\Cal D}
\define\E{\Cal E}
\define\U{{\Cal U}} 
\redefine\c{$\clubsuit$} 
\def\<{\langle}
\def\>{\rangle}
\def\chip{\chi^\prime}
\def\cfp{{\Cal CF}^\p}
\def\hcfp{{\Cal HCF}^\p}  
\def\cfpx{{\Cal CF}^\p_\xi}
\def\hcfpx{{\Cal HCF}^\p_\xi}

\define\pd#1#2{\dfrac{\partial#1}{\partial#2}}

\nopagenumbers

\define\Ca{1}
\define\FF{{\Cal F}}
\define\LL{{\Cal L}}
\define\MM{{\Cal M}}
\define\WW{{\Cal W}}
\define\TT{{\Cal T}}

\redefine\o{{\omega}}
\define\oa{\omega_{\alpha}}
\redefine\a{{\alpha}}
%\redefine\c{{\circ}}
\define\Mw{M_{\omega}}
\define\Sw{\Sigma_{\omega}}
\define\Gw{\Gamma_{\omega}}
\define\Ge{\Gamma_{\eta}}
\define\Gf{\Gamma_{\bar f}}

\define\Gwb{\bar\Gamma_{\omega}}

\define\codim{\operatorname{codim}}
\define\Ha{3}
\define\Bo{2}

\heading  Statement of results \endheading

\subheading{1. Introduction} The Morse theory of closed 
$1$-forms was begun by S. P. Novikov \cite{N}, \cite{N1}; 
a survey can be found in \cite{P}. The 
standard Morse inequalities for functions were generalized, by 
Novikov, to 
forms using a twisted cohomology defined by the $1$-form in place 
of the usual cohomology. In some cases it is known that these 
inequalities are sharp \cite{F}, \cite{P}. 

In this paper we will address the problem of whether it is 
possible to improve these Novikov inequalities (or any other 
inequalities
of a similar nature) {\it assuming additionally
that the given
closed 1-form is harmonic with respect to some Riemannian metric.} 
This problem was first considered in the paper of E.Calabi \cite{C}
even before  Novikov's theory. 
There is one obvious restriction -- namely,
any non-constant harmonic form has no local minima and maxima -- 
and Calabi
asked in \cite{C} if there are any  
other restrictions on the critical point structure of Morse 
harmonic $1$-forms.  

Motivated by this , Calabi
solved instead in \cite{C} another fundamental problem: he gave a 
complete 
topological criterion for a closed 
Morse 1-form to be harmonic with respect to some Riemannian 
metric.
We will call such forms {\it intrinsically harmonic } or {\it Calabi 
forms}.
We will describe the theorem of Calabi briefly in Section 3.

The thrust of the main results in this paper is that, at least under 
certain quite general conditions,
{\it there exist no  
restrictions on the Morse numbers of
harmonic forms except 
for the trivial restrictions mentioned above}. 

The conditions we need to
impose are closely tied to the topological nature of the (singular) 
foliation which is associated to a Morse $1$-form (see below for 
the definition). In particular, compactness of the leaves is a crucial
item. Our results will show that if the foliation defined by a form 
$\o$
 is made up entirely of {\it non-compact} leaves or, at the other 
extreme,
 belongs to a cohomology class which is a scalar multiple of an 
integral
class (in which case all the leaves of the foliation are compact),
then there is a Calabi form $\tilde\o$ which is {\it contiguous} 
to  
$\omega$, i.e. it is in the same cohomology class
 and has the
 same number of critical points of each index.

We will also present some results concerning two issues related to 
the
 existence of compact leaves in the foliation of a closed Morse $1$-
form.
 In section 9 we show that such a form defines a natural 
decomposition of the 
 manifold into pieces, each of whose interiors consists of all 
compact or all
 non-compact leaves. The combinatorics of this decomposition plays 
a crucial
 role in determining whether the form is contiguous to a Calabi form.

In Sections 10-11 we present two results which give restrictions on 
the
 cohomology class of a form when either all or none of the leaves of 
its
 foliation are compact. In Section 12 we describe several examples 
for surfaces
 which will illuminate and illustrate our results.

It is a theorem of Lalonde-Polterovich~\cite{LP} that, given any 
non-trivial class 
in $H^1(M,\R)$, then, for a generic Riemannian metric, the harmonic
form representing that class has only Morse singularities. 
In other words, generically harmonic forms have only Morse type 
singularities. This allows us to avoid considering more complicated
singularities and leaves with the question addressed here of how 
many critical points these Morse harmonic forms can have. 

Finally we would like to thank Robert Kotiuga for bringing the work of
Calabi to our attention and for many stimulating discussions.

\subheading{2. The main results} Here we formulate in a precise 
form
our principal results.

\proclaim{ Theorem 1} Suppose that $\o$ is a closed Morse $1$-form 
on a closed manifold $M$ which represents a scalar multiple of an 
integral
cohomology class and such that $\o$ has no critical points
of index $0$ or 
$n=\dim M$. Then there exists an intrinsically harmonic Morse 
1-form $\tilde \o$ 
which has the same numbers of critical points of all indices as $\o$ 
and which
is in the same cohomology class as $\o$. \endproclaim

This theorem implies that, if there would exist Morse type
inequalities for harmonic forms estimating the numbers of critical 
points by
information depending only on the cohomology class of the form,
they have to be also true for any closed 1-form representing a scalar
 multiple of
an {\it integral} cohomology class. We refer to such cohomology 
classes as
 being of {\it 
rank one}. 
Following S.Novikov~\cite{N}, we say that a cohomology class 
$\xi\in H^1(M,\R)$ 
{\it is of rank $k$} if the image of the homomorphism $\pi_1(M)\to 
\R$ 
determined by $\xi$ is a free abelian group of rank $k$.

We recall from \cite{N} that, for a closed 1-form $\o$ of rank one 
on a closed manifold $M$, the Novikov inequalities take the form:
$$m_p (\o )\geq b_p ([\o ])+q_p ([\o ])+q_{p-1}([\o ]), $$
where $m_p (\o )$ is the number of critical points of index $p$ and 
$ b_p ([\o ]), q_p ([\o ])$ are the rank and torsion Novikov 
numbers of the
cohomology of $M$ with local coefficients defined by $[\o ]$ (see 
\cite{N} 
and \cite{F} for precise definitions). In \cite{F} it is shown that 
these
inequalities are exact. Combining this with Theorem 1 gives:

\proclaim{Corollary} If $\pi_1(M)=\Z$ and $n=\dim M \ge 6$, then 
for 
any nonzero class $\xi\in H^1 (M;\R )$, there 
exists an intrinsically harmonic Morse form $\o$ representing $\xi$ 
such
that $m_p (\o )=b_p (\xi)+q_p (\xi )+q_{p-1}(\xi )$.

\endproclaim
This follows by combining the main Theorem of \cite{F} and Theorem 
1 above
and observing that, since $b_0(\xi)=q_0(\xi)=b_n(\xi)=q_{n-
1}(\xi)=0$, 
the 1-form
given by \cite{F} 
realizing the given non-zero cohomology class $\xi$ has no critical 
points of
indices 0 and $n$.

It is an easy fact
that rank one Morse forms must have all compact leaves.
 In the next theorem we consider the opposite situation: we
assume that all leaves are non-compact. In addition we need to 
assume that $\o$ is {\it generic}, which we define to
 mean that any singular leaf contains only one critical point. It is 
easy
 to see that any Morse form can be perturbed slightly so that it 
becomes generic
 without changing the numbers of critical points or the cohomology 
class.

\proclaim{Theorem 2} If $\o$ is a {\it generic} closed Morse $1$-
form all of 
whose leaves are non-compact then $\o$ is intrinsically harmonic.
\endproclaim

Theorem 2 is false if $\o$ is not generic. See example (4) in 
Section
12.

Theorems 1 and 2 will be proved in Sections 7 and 8.

\heading   Calabi's theorem     \endheading  

\subheading{3} In this section we will recall Calabi's Theorem 
\cite{C}.

Let $\omega$ be a fixed closed 1-form on a closed manifold $M$.

A smooth path $\gamma:[0,1]\to M$ will be called {\it $\omega$-
positive} if 
$\omega(\dot\gamma(t))>0$ for any $t\in [0,1]$. 
Calabi considered the following condition on $\omega$ which he 
calls
{\it transitivity}:

{ (a)}  For any ordered pair of points 
 $x$ and $y$ in $M$, which are not singular points of
$\omega$, there exists an $\omega$-positive path from $x$ to $y$.

Evidently, (a) implies: 

{ (b)}  For any non-singular point $x \in M$ there 
exists a {\it closed $\omega$-positive} path $\gamma$ through $x$.

We shall see that, in fact, (a) is equivalent to (b). 

The main result of \cite{C} can be reformulated as follows. 

\proclaim{Theorem (\cite{C})} A closed Morse 1-form 
$\omega$ on $M$ is harmonic with respect to some Riemannian 
metric
if and only if condition (b) holds. 
\endproclaim
 
\subheading{Remark} It is pointed out in \cite{C} 
 that (a) precludes the existence of 
critical points of index $0$ or $n=\dim M$ and that non-singular 
forms satisfy (a).

\heading  Singular foliation of a Morse form \endheading

{\subheading 4} Let again $\omega$ denote a closed 1-form with 
only Morse
type singularities. 

For any simply-connected open set $U\subset M$ we have 
$\omega|_U=\,df_U$
for some Morse function $f_U:U\to \R$ determined up to a constant. 
We will consider the foliation in $U$ determined by the level sets of 
$f_U$.
Choosing a covering $\{U\}$ of $M$ we see that all these foliations
match together to form a foliation $\FF$ of $M$. We will call it 
{\it the foliation determined by the form $\omega$.}

Note this foliation has finitely many singular points which are 
the critical points of $\omega$. Locally the structure of the singular
foliation $\FF$ around a critical point of index $d$ has the form
$$-x_1^2-\dots -x_d^2+x_{d+1}^2\dots +x_n^2\ =\ c.$$

A {\it leaf of $\FF$} is, by definition, any maximal subset $\LL$ of 
$M$
such that for any two points $x,y$ in $\LL$ there exists a smooth 
path
$\gamma:[0,1]\to M$ which connects $x$ with $y$ and such that 
$\omega(\dot\gamma(t))=0$ for all $t$. There are finitely many 
leaves
containing the singular points; we will call those leaves {\it 
singular}.

 Note that 
the {\it non-singular} leaves, i.e. those 
containing none of the  critical points of $\o$, are smooth while the 
{\it singular} leaves are smooth except at each of the critical points 
which has a neighborhood homeomorphic to a cone over $S^{d-
1}\times S^{n-d-1}$, where $d$ is the index of the critical point. 
Each nonsingular leaf $\LL$ has a natural topology 
(which will be called {\it the leaf topology}) and smooth 
manifold structure such that the inclusion $\LL\to M$ is an 
immersion;
similarly for the singular leaves. We 
will be particularly interested in those leaves containing critical 
points of index $1$ or $n-1$. In these cases removing a critical point 
locally disconnects the leaf and we will find it useful to define a 
{\it singular leaf component} to be the closure (in the leaf topology)  
of a connected component of $\LL-\Sigma(\LL)$, where 
$\Sigma(\LL)$ 
denotes the set of all singular points of the leaf $\LL$.

\heading    Calabi graphs     \endheading

{\subheading 5} We now make some preliminary constructions to set 
up the 
proof of Theorem 1. In the next section 6 we will see that the Calabi
properties of closed 1-forms can be interpreted by analogous 
properties
of oriented graphs.

Let $\Gamma$ be an {\it oriented} connected finite 
graph. Consider the following two properties. 

\item{ (a')}  If $x, y\in\Gamma$ (it suffices to consider only 
vertices) 
then there is a path from $x$ to $y$ that traverses edges of 
$\Gamma$ only in the positive direction. 

\item{ (b')}  For any point $x\in\Gamma$, there exists a 
closed path through $x$ 
that traverses edges of $\Gamma$ only in the positive 
direction. 

Clearly (a') implies (b').

\proclaim{ Lemma 1} If $\Gamma$ satisfies (b'), then it satisfies 
(a').
\endproclaim
\demo{Proof} Condition (b') tells us that $\Gamma$ is the union of 
closed edge-paths which can be traversed in a positive direction 
with respect to the orientation. Since $\Gamma$ is connected, then, 
for any two points $x, y\in\Gamma$, there is a sequence of points 
$x=x_1, x_2,\cdots, x_k =y$ such that $x_i$ and $x_{i+1}$ are on a 
common closed edge path $C_i$. Clearly $x_i$ and $x_{i+1}$ can be 
joined by a positive path on $C_i$. Putting these together gives the 
desired path from $x$ to $y$. \qed
\enddemo

\proclaim{Definition}  We shall say that $\Gamma$ is a 
{\it Calabi graph} if it satisfies (a')\ $\sim$ (b').
\endproclaim

In Figure 1 we show an example of a Calabi graph and a non-Calabi 
graph.

\qquad\qquad\quad

\centinsert{\pscaption{\psboxscaled{500}{figT}}{Figure 1}}
%\botcaption{FIGURE 1}
%\endcaption
\bigskip
\subheading{6} With any closed Morse 1-form $\omega$ on a closed 
manifold $M$, 
generating a singular foliation with {\it all compact leaves}, we 
shall 
associate an oriented graph $\Gamma_{\omega}$. With this goal in 
mind, 
we introduce an equivalence relation in $M$. We declare two 
points in $M$ to be equivalent if they lie on the same leaf of the 
foliation 
$\FF_{\omega}$, generated by $\omega$. When all the leaves of 
$\FF_{\omega}$ are compact, then the 
singular leaves are isolated (i.e. each singular leaf has a 
neighborhood, free of other singular leaves) and 
the quotient space $\Gamma_{\omega} := M/\sim$ 
is a graph. Indeed, any non-singular leaf $\LL$  has a collar, 
consisting of non-singular leaves. Hence, its neighborhood in 
the quotient space $\Gamma_{\omega}$ is an open interval. 
Suppose that $\LL_0$ is a singular leaf corresponding to a point 
$v\in\Gw$. The form $\o$ is exact in some neighborhood of $\LL_0$. 
Suppose $\o =\,df$ and $f^{-1}(a)=\LL_0$. Then, for some $\epsilon$, 
$f|_{f^{-1}(a,a+\epsilon)}$ and $f|_{f^{-1}(a-\epsilon ,a)}$ 
are fibrations and the fibers are non-singular leaves. Thus the 
complement of $\LL_0$ in this neighborhood is a product of an open 
interval with a finite number of compact manifolds which define a 
finite number of edges of $\Gw$ with $v$ viewed as a vertex.

Suppose $x_i \in\LL_0$ is a critical 
point of index $s_i$. If $1<s_i <n-1$, then it follows from the local 
structure of a non-degenerate critical point that the intersection 
of each leaf with a small neighborhood of $x_i$ is connected. We 
conclude 
that if all the indices $s_i$ for a given singular leaf $\LL_0$ 
satisfy $1<s_i <n-1$, then the leaves on both sides of $\LL_0$ are 
connected and so $v$ 
is just an interior point of an edge of $\Gw$. If $s_i =1$ or $n-1$ 
 and $n>2$, 
then locally we have two non-singular leaves coalescing into the 
singular leaf from one side and one emerging on the other side (see 
Figure 2). (Of course it is possible that the what looks like two 
leaves locally may actually be part of the same leaf globally.) 

\qquad\qquad
\centinsert{\pscaption{\psboxscaled{500}{fig2}}{Figure 2}}
\bigskip

Thus when one or more of the indices $s_i$ equal 1 or $n-1$, then 
$v$ may 
be a {\it true} 
vertex of $\Gw$. 
If $s_i =0$ or $n$, then we have 
non-singular fibers on only one side and $v$ will again be a {\it true} 
vertex.

In Figure 3 we show two examples of Morse $1$-forms with all 
compact leaves whose associated graphs are those in Figure 1. The 
forms are the pullbacks of the canonical $1$-form $\,d\theta$ on 
the 
circle $S$ via the obvious projection maps onto $S$.
\qquad\qquad\quad
\centinsert{\pscaption{\psboxscaled{500}{fig3}}{Figure 3}}
\bigskip
The following lemma is obvious from the preceding discussion.

\proclaim {Lemma 2} $(M, \omega)$ (with all leaves compact) 
 satisfies properties 
(a) and (b) of section 3 iff the associated oriented graph $\Gw$
 satisfies properties (a') and (b') of section 5, respectively. 
\endproclaim

\heading  Proof of theorem 1 \endheading

\subheading{7}
First, we remark that we may assume that $n=\dim M > 2$ without 
loss
of generality. In fact, if $n=2$, then from the Euler-Poincare 
Theorem
it follows that any 1-form which has no maxima and minima, has 
precisely
$-\chi(M)$ critical points of index 1. Thus we will assume in the 
sequel
that $n>2$. In fact, with some extra care, all our arguments will 
work in dimension $2$. Our figures are all in dimension $2$ so the 
reader should interpret them accordingly.  
 
Multiplying the form $\omega$ by a scalar, we may assume that
it is the differential of a (multivalued) Morse function $f:M\to S$, 
where $S$ 
is an oriented circle. 
More precisely, this means that $\o\ =\ f^\ast(d\phi)$, where 
$d\phi$ is the 
standard angular form on the circle.

We may also assume, after a perturbation, that all the 
critical values of $f$ are distinct. Then, since there are no critical 
points of index $0$ or $n$, every vertex of the graph $\Gw$ (cf. 
section 6) is 
trivalent.

Consider the map $\psi_f :\Gw\to S$ induced by $f$. 
Define the {\it complexity} of $\o$ to be the smallest cardinality of 
any $\psi_f^{-1}(a)$, where $a$ ranges over the regular values of 
$f$. For example in Figure 3   
the complexities of the two graphs are 
$2$ and $1$, respectively. 

We will show that either:
\roster
\item $\Gw$ has no vertices, or
\item there exists another closed $1$-form $\o'$ with strictly 
smaller complexity such that $\o$ and $\o'$ are {\it 
contiguous}.
 Recall, that
this means, that $\o$ and $\o'$
have the same numbers of critical
points of all indices and belong to the same cohomology class.
\endroster
Thus, using induction , we will eventually reach case (1)
i.e. we will find a closed 1-form $\tilde \o$ contiguous to $\o$ and 
such  
that $\Gamma_{\tilde\o}$ is a circle 
and $\psi_f$ is a covering map. In this case it is clear that 
$\Gamma_{\tilde\o}$ is 
a Calabi graph and we are done.

Suppose that $\Gw$ has at least one vertex. Choose a regular value 
$a$ 
such that the cardinality of $\psi^{-1}(a)$ equals the complexity of
$\o$. Cut open $M$ 
along $M_a =f^{-1}(a)$ to obtain a manifold $\bar M$. The map $f$ 
induces 
a Morse function $\bar f:\bar M\to [0,2\pi]$ such that the two 
identical 
boundary components $M_a^+$ and $M_a^-$ map to $0$ and $2\pi$ 
respectively. The leaves of the foliation defined by $\bar f$ are the 
same as those of $\o$ and the graph $\Gf$ it defines (cf. section 6) 
can be obtained 
from $\Gw$ by cutting some of the edges.
We now change the function $\bar f$, using the Morse-Smale 
reindexing technique (see \cite{M}), to obtain a self-indexing Morse 
function $\bar g:\bar M\to [0,2\pi]$ with the same number of 
critical 
points of each index as $\bar f$ and such that $\bar g$ agrees with 
$\bar f$ in a neighborhood of $\partial\bar M$. The self-indexing 
property means that for any two critical points $x, y$, 
index$(x)>$index$(y)$ implies $\bar g(x)>\bar g(y)$. We may also 
assume that all the critical values of $\bar g$ are distinct. If we 
reglue $M_a^+$ and $M_a^-$, the function $\bar g$ induces a new 
Morse 
$1$-form $\o'$ on $M$ which is contiguous to $\o$. See Figure 4.  

%\qquad\qquad\qquad\qquad
\centinsert{\pscaption{\psboxscaled{500}{figC}}{Figure 4}}
\bigskip

Now consider the associated 
oriented trivalent graphs $\Gamma_{\o'}$ and $\Gamma_{\bar g}$. 
The vertices of $\Gamma_{\bar g}$ are of two types: those with two 
ingoing and one outgoing edge (they necessarily
belong to singular leaves containing
critical points of index 1) and those with one ingoing and two 
outgoing edges (they necessarily belong to singular leaves 
containing 
critical points of index $n-1$). The self-indexing property implies 
(since we assume that $n>2$ and so $1<n-1$ 
that 
the induced 
map $\psi_{\bar g}:\Gamma_{\bar g}\to [0,2\pi]$ has the property 
that 
there is some $b\in (0,2\pi)$ such that $\psi_{\bar g}$ maps each 
vertex of
the first type into $(0,b)$ and each vertex of the second type into 
$(b,2\pi)$. 
As a consequence, if there are any vertices at all, then the 
cardinality of 
$\psi_{\bar g}^{-1}(b)$ is strictly smaller than that of $\psi^{-
1}(a)$. 
But this means that the complexity of $\o'$ is smaller than 
the complexity of $\o$. 

For example, applying this modification to 
the non-Calabi example in Figure 3 and its graph in Figure 1 will 
produce the form given by the other example in Figure 3.

This completes the proof. \qed

\heading Proof of Theorem 2\endheading

{\bf 8.}  Suppose the Morse $1$-form 
$\o$ has all non-compact leaves. We will show that it satisfies the 
Calabi condition (b) of section 3. Let $x\in M$. If $x$ belongs to a 
non-singular 
leaf, we will
denote it by $\LL_x$; if $x$ belongs to a singular leaf, then we will 
denote
by $\LL_x$ the {\it singular leaf component} (cf. section 4) 
containing $x$.

By statement (2) of Proposition 1 (cf. below in section 9), since $\o$ 
is assumed
 to be 
generic, we know that $\LL_x$ is non-compact. 

First of all we point out that if any non-singular  
point $x'$ of $\LL_x$ satisfies Calabi condition (b), i.e. there is a 
$\o$-positive closed path $\gamma$ passing through $x'$, then any 
other 
non-singular point $x''$ in $\LL_x$ satisfies Calabi condition (b) as 
well. 
To show this we
connect $x''$ to $x'$ by 
a path in $\LL_x$. Now, $\gamma$ crosses $\LL_x$ transversely at 
$x'$ and so we can pull $\gamma$, as it passes through $x'$, along 
$\gamma$, so that it now passes through $x''$. It is clear that we 
may keep $\gamma$ positive during this deformation. Thus, it 
suffices to find an $\o$-positive closed path intersecting $\LL_x$ 
anywhere. 

Since $\LL_x$ is not compact there exists a limit point $y$ of 
$\LL_x$ not in $\LL_x$. Since we have explicit models for the 
foliation near $y$, for the singular and non-singular case, we can 
see that all points on $\LL_y$ close to $y$ are also limit points of 
$\LL_x$.  So we may assume that $y$ is not a critical point. If 
$U$ 
is a sufficiently small coordinate neigborhood of $y$, then the 
foliation on $U$ consists of parallel hyperplanes and so $\LL_x \cap 
U$ must contain a sequence of these hyperplanes converging to the 
one containing $y$. Choose any two points $x', x''$ of $\LL_x 
\cap U$ in different hyperplanes. We can connect these two points by 
a path $\rho:[-1,1]\to \LL_x$ with $\rho(-1)=x''$ and $\rho(1)=x'$.
We can also obviously connect them 
by a positive path $\xi$ in $U$ (say, running from $x'$ to $x''$). We 
now push $\rho$ slightly off $\LL_x$ to make it into a positive 
path. To do so choose a vector field $X$ on $\LL_x$ so that $\o (X)$ 
is 
a positive constant. We can assume $X$ is tangent to $\xi$ at 
$x'$ and $x''$.  Define 
$\eta (t)=\gamma (t)+\epsilon tX(\gamma (t))$ for $\epsilon$ small 
enough, so that $\eta (-1)$ lies ahead of $\eta (1)$ on $\xi$. Now we 
can put $\eta$ together with the portion of $\xi$ that runs from 
$\eta (1)$ to $\eta (-1)$ to create a positive closed path which 
intersects $\LL_x$ at $\eta (0)$. Of course, we need to round 
corners. 
See Figure 5.  \qed

\qquad\qquad\qquad
\centinsert{\pscaption{\psboxscaled{500}{fig4}}{Figure 5}}
\bigskip

\subheading{9. Compact and non-compact leaves} 
It is clear from Theorems 1 and 2 that the general problem of 
deciding
 whether a given closed Morse 1-form is contiguous to a Calabi form 
may 
depend
on the structure of the foliation defined by $\o$. We will now
 explain how the manifold $M$ breaks up, in a nice way, into pieces 
made
 up (roughly) 
of all compact leaves or all non-compact leaves. 

It is a theorem of S. Novikov, mentioned in \cite{T}, that when  $\o$ 
is 
non-singular (and $M$ is connected) then either all the leaves are
compact, or all the leaves are non-compact. This is certainly not 
true for general foliations (e.g. the Reeb foliation) but, in 
Proposition 1, we will generalize Novikov's theorem to foliations 
coming from a closed Morse $1$-form, showing that the compact and 
non-compact leaves give a nice decomposition of $M$.   
The terms used in the statement were defined above in section 4 and 
in 
section 2.

\proclaim{Proposition 1}    If $\o$ is a closed Morse $1$-form on 
$M$, 
then
$M$ is the union of two compact $n$-dimensional submanifolds 
$M_c$ and
 $M_{\infty}$, with a common (but generally singular) boundary, 
satisfying:
\roster
\item $\Int M_c$ is a union of all of the compact leaves plus some 
compact 
singular leaf components of non-compact leaves;
\item $\Int M_{\infty}$ is a union of all non-compact non-singular 
leaves , all non-compact singular leaf components and 
some compact singular leaf components of non-compact leaves.
\item $\partial M_c =\partial M_{\infty}=M_c \cap M_{\infty}$ is a 
subvariety which is smooth except at a finite number of points, 
which are critical points of $\o$, and is a union of some compact 
singular 
leaf components of non-compact singular leaves. 
\item If $\o$ is generic, then
$M_c \cap M_{\infty}$ contains all compact singular leaf 
components of non-compact
leaves.
\endroster
\endproclaim

The theorem of Novikov follows from Proposition 1 since if there 
are
no critical points, then by statement (3) 
$$\partial M_c\ =\ \partial M_\infty \ =\ \emptyset$$
and so $M=M_c$ or $M=M_\infty$.

\demo{Proof}    
(1) Let $U$ be the union of all the compact leaves of the foliation 
determined by the form $\o$. We show that $U$ is open. If $\LL$ is a 
compact leaf, then it has a neighborhood $W$ with compact closure 
and containing no critical points other than those on $\LL$, such that 
$\o |\bar V=\d f$ for some smooth function on $\bar V$ with $f^{-
1}(0)=C$. Then there exists $\epsilon >0$ such that $\overline{f^{-
1}(-\epsilon ,\epsilon )}\subseteq V$ and so $f^{-1}(a)$ is a compact 
leaf for $|a|<\epsilon$.

(2) Let $X$ be the union of all compact leaves and compact leaf 
components of non-compact leaves. Then $X$ is closed. Indeed if we 
consider the non-singular foliation $\FF_0$ defined by $\o$ on $M_0 
=M-$critical points of $\o$, we can apply a theorem of 
Haefliger~\cite{H, p.386} to conclude that the union $X_0$ of all the 
closed leaves of $\FF_0$ is a closed subset of $M_0$. But then $X$ is 
just the closure of $X_0$ in $M$.

Let $M_c =\bar U$ and $V=M-\bar U$. From (2) we conclude that 
$Y=M_c -U$ is a union of some compact leaf components of 
non-compact leaves. $Y$ is a codimension one submanifold of $M$ except 
for singularities at the critical points of $\o$ in $Y$. If $x$ is such a 
critical point of index $q$ then, for some neighborhood $W$ of $x$, 
$(W,Y\cap W)\cong\text{ cone over }(S^{n-1},S^{q-1}\times S^{n-q-
1})$ if $1<q<n-1$, while if $q=1\text{ or }n-1$ then $(W,Y\cap 
W)\cong\text{ cone over }(S^{n-1}, S^{n-2})\text{ or }(W,Y\cap 
W)\cong\text{ cone over }(S^{n-1},S^0\times S^{n-2})$. In this last 
case we call $x$ {\it special }. It is clear that $M_c$ is a manifold 
with boundary $Y$ which is smooth except at the critical points and 
at those points we have a very explicit model of the singularity. At 
every point of $Y$, except the special points, $Y$ locally separates 
$M$ into two components exactly one of which must lie in $U$ and 
the other lies in $V$. At the special points $Y$ locally separates $M$ 
into three components if $n>2$ or four components if $n=2$. If $n>2$ 
then either one or two of the components lie in $U$. See Figure 6.

\centinsert{\pscaption{\psboxscaled{500}{figA}}{Figure 6}}
\bigskip

It remains to discuss the generic case. Suppose $C$ is a compact 
leaf component of a non-compact leaf $\LL$. Then $C$ contains 
exactly one critical point $x$ which must be of index $1$ or $n-1$. 
As above we have a neighborhood $W$ of $x$ so that $(W,\LL\cap 
W)\cong\text{ cone over }(S^{n-1},S^0\times S^{n-2})$ where $C\cap 
Wgong\text{ one of the two cones over }S^{n-2}$. Now we can 
choose a regular neighborhood $R$ of $C$ so that $R=W$ near $x$ but 
$(R-W,C-W)\cong((C-W)\times [-1,1], (C-W)\times 0)$ since there 
are no other critical points on $C$. See Figure 7. It is clear that the 
leaves of the foliation defined by $\o$ on one side of $C$ must lie 
entirely in $R$ and are homeomorphic to $C$. Since $C$ is compact 
this show $C\subseteq\bar U$. \qed

\centinsert{\pscaption{\psboxscaled{500}{figB}}{Figure 7}}
\bigskip

\enddemo

{\bf 10.} The next two propositions describe some relations between 
the
 cohomology class
represented by a Morse 1-form and the presence of compact leaves 
in its
 foliation.

\proclaim{Proposition 2} A class $\xi\in H^1 (M;\R )$ has a 
representative Morse $1$-form with all compact leaves if and only 
if the homomorphism
$\pi_1(M)\to \R$ determined by the cohomology class $\xi$ can be 
factorized 
$\pi_1(M)\to F\to \R$
through a free group $F$.
\endproclaim 

It will be convenient for us to refer to the cohomology classes 
satisfying
the condition of Proposition 2 as {\it split} cohomology classes. A 
more 
geometric description of this condition is the following. The 
manifold $M$
can be cut into several pieces by disjoint codimension one 
submanifolds
so that the restriction  of the class $\xi$ to each piece is of rank 
one.

\remark{ Remarks} 
\item{(i)} As pointed out in Section 2, a rank one class has the 
property that
 {\it every} representative Morse 1-form has all compact leaves. In 
Section
 12 we will give examples that suggest that only rank one classes
 have this stronger property.
\item{(ii)}  
It seems to be difficult to give a purely cohomological 
criterion for $\xi$ to be split, but a necessary condition is that 
there exist a finite subset $\{\xi_i \}
\subseteq H^1 (M;\Z )$ such that $\xi$ is a linear combination of the 
$\xi_i$ 
and all Massey products (see \cite{D}) of the $\{\xi_i \}$ vanish. 
However {\it this condition is not sufficient}. An example can be 
obtained as 
follows. Choose a slice link $L$ of two 
components which is not a 
homology boundary link. 
According to \cite{Hi, p.72}    the link in Figure 8 
is such an example. Define $M$ to be the result 
of $0$-surgery on $S^3$ along $L$. Then $H^1 (M;\Z )\cong\Z^2$ and 
all Massey products vanish. But there is no map $\pi_1 (M)\to F,\ F$ 
a 
free group inducing an injection on $H^1$ and so, any class $\xi$ 
of rank two is not split.

\endremark

\qquad\qquad\qquad\qquad
\centinsert{\pscaption{\psboxscaled{500}{fig5}}{Figure 8}}
\bigskip

\demo{Proof of Proposition 2}  Let $\o$ be a Morse $1$-form 
representing $\xi$ with 
all compact leaves and $\Gw$ the associated oriented graph as 
defined in section 6.
Let 
$\phi~:M\to\Gw$ be the identification map. It is clear that 
$\xi\in\phi^{\ast}H^1 (\Gw ;\R )$ and half of Proposition 2 follows
since any connected graph has a free fundamental group.

Now suppose we are given $\xi$ and a map $f:M\to W$, where $W$ is 
a 
wedge of circles, such that 
$\xi=f^{\ast}(\xi ')$ for some $\xi'\in H^1 (W;\R )$. 
Choose a point $a_i$ in each 
circle $S_i$ of $W$; we may assume that $a_i$ is a regular value of 
$f$. Then $M_i =f^{-1}(a_i )$ is a smooth compact submanifold of 
$M$. Let $\bar M$ be the result of cutting $M$ open along the $M_i$. 
Then $\partial\bar M$ consists of two copies $M_i^+$ and $M_i^-$ of 
each $M_i$. Obviously $\xi |_{\bar M}=0$. Now assume each $S_i$ is 
oriented and set $p_i =\<\xi ',[S_i]\>$, the period of the class $\xi'$ 
on 
the circle $S_i$. We may assume that the 
orientations are such that a short path in $M$ which crosses $M_i$, 
hitting $M_i^+$ before $M_i^-$, maps into $S_i$ with positive 
orientation. Now choose a Morse function $g:\bar M\to\R$ with no 
critical points on $\partial\bar M$, which is constant on each 
boundary component so that $g(M_i^+ )-g(M_i^- )=p_i$. We may also 
assume that $g|_{\partial\bar M}$ fits together smoothly to form a 
smooth $1$-form $\o$ on $M$. Clearly $\o$ is Morse and closed, has 
zero periods in $\bar M$ and has period $p_i$ around any closed 
curve which crosses $M_i$ \lq\lq positively'' but no other $M_j$. 
Thus 
$\o$ represents $\xi$. Furthermore, the leaves of the foliation 
defined by
$\o$ are components of $f^{-1}(a)$ for some $a\in\R$ and so they are 
compact. 
\qed
\enddemo

{\bf 11.}
\proclaim{Proposition 3} If $\o$ is a closed intrinsically harmonic
Morse $1$-form 
whose associated foliation has a compact leaf, then there exists 
some {\it non-zero } $\theta\in H^1 (M;\Z )$ such that $\theta\cup 
[\o ]=0  \in H^2(M,\R)$.
\endproclaim

We will show in Section 12 that Proposition 3 is false if $\o$ is not
 assumed to be intrinsically
harmonic. 

\demo{Proof} 
If $\o$ has a compact leaf, then, by statement (1) of Proposition 1,
there exist {\it non-singular} compact leaves. 
Suppose $\LL$ is a non-singular compact leaf of $\o$. Then 
$\LL$ is a closed submanifold of $M$ and the normal bundle to $\LL$
is trivialized by the form $\o$. Also, $\LL$ cannot bound in $M$ 
since
this would contradict the Calabi condition. 

We may find a cylinder $U = [-1/2, 1/2] \times \LL$ in $M$ fibred by
the leaves of the foliation determined by $\o$, such that the 
initially
choosen leaf $\LL$ is identified with $0\times \LL$. Now, let us 
construct a closed
1-form $\alpha$ on $M$ such that:
   
(1) $\alpha$ vanishes outside the cylinder $U$;

(2) on the cylinder $U$ the form $\alpha$ is given by
$\alpha|_U = d(\phi\circ p)$ where $p:U\to  [-1/2, 1/2]$ denotes the 
projection,
and $\phi:  [-1/2, 1/2]\to  [-1/2, 1/2]$ is a smooth function 
satisfying:

$$\phi(t)=t \quad\text{for all}\quad |t|< 1/2 - \epsilon$$
$$ \phi(t)=1/2 \quad\text{for all}\quad t\in [1/2-\delta, 1/2]$$
$$\phi(t)= -1/2\quad\text{for all} \quad t\in [-1/2, -1/2+\delta]$$
where $0<\delta<\epsilon$. 

Thus we obtain that $\alpha$ represents a nonzero integral 
cohomology
class $\theta\in H^1(M,\Z)$ and $\alpha\wedge \o = 0$. Therefore, on 
the level
of cohomology, we get $\theta\cup [\o] = 0$ in $H^2(M,\R)$.\qed
\enddemo

Suppose, for example, that the pairing 
$H^1 (M,\R)\times H^1 (M,\R)\to H^2 (M,\R)$ is non-degenerate, i.e. 
for any
$\xi\in H^1(M,\R)$ there exists $\eta\in H^1(M,\R)$ such that 
$\xi\cup\eta\ne 0$.
Then it 
follows
from Proposition 3 that, if a class $\xi\in H^1(M,\R)$ is completely 
irrational, then any 
representative 
harmonic Morse $1$-form has all non-compact leaves.
(We  say that a cohomology class $\xi$ is {\it completely 
irrational} 
if its rank is equal to the first Betti number of $M$.)
This applies to the cases when 
$M$ is 
an orientable surface or $M=S^1 \times\cdots\times S^1$.

\heading examples\endheading

{\bf 12.} We complement and illustrate Theorem 2 and Propositions 
2 and 3
with some examples in dimension two. We first summarize the
assertions of these results for a class $\alpha\in H^1 (M;\R )$, 
where
$M$ is a closed orientable surface, and a representative Morse
$1$-form $\oa$.
\roster
\item If rank $\alpha =1$, then any $\oa$ has all compact leaves.
\item $\alpha$ is split if and only if there exists some $\oa$ with 
all compact
 leaves.
\item If $\alpha$ is completely irrational and $\oa$ is Calabi, then 
no leaves are compact.
\item If $\oa$ is generic and has no compact leaves, then $\oa$ is 
Calabi.
\endroster

   We will construct the following four examples of $\alpha$ and 
$\oa$ which will illustrate the sharpness of these results.
\proclaim{Examples}
\item{(1)} Rank $\alpha=$ one less than maximal and $\oa$ Calabi 
with some
compact leaves.
\item{(2)} Rank $\alpha >1, \alpha$ split and $\oa$ Calabi with all 
leaves
 non-compact.
\item{(3)} $\alpha$ completely irrational and $\oa$ non-Calabi with 
some compact leaves.
\item{(4)} $\oa$ non-Calabi but all leaves non-compact (and so 
$\oa$ is not generic).  
\endproclaim

We first describe three general connected sum operations which 
accept as input
 two closed $2$-dimensional manifolds $M_1 ,M_2$ with closed Morse $1$-forms
 $\o_1 ,\o_2$ and constructs a closed Morse $1$-form $\o$ on the
 connected sum $M_1\sharp M_2$ such that $\o |M_i -D_i =\o_i |M_i 
-D_i$,
 where $D_i$ is any disk in $M_i$ containing no critical points of
 $\o_i$ and within which the connected sum will be constructed.   
 
Choose smooth functions $f_i :D_i\to\R$ such that $\o_i =d f_i$. 
Suppose that
we can identify each $D_i$ with an open rectangle $(a_i ,b_i )\times
 (c_i ,d_i )$ in $\R^2$ so that $f_i (x,y)=y$. Note that we can then
 move the rectangle by any translation in $\R^2$.
 Now we can perform a connected sum of
 $D_1$ and $D_2$ ambiently by connecting them with a straight tube 
in
 $\R^3$ in three different ways.
\item{A.} $f_1 (D_1 )\cap f_2 (D_2 )=J$ is non-empty, the 
connecting
 tube intersects $D_i$ inside $f_i^{-1}(J)$ and is approximately
horizontal.
\item{B.} $f_1 (D_1 )\cap f_2 (D_2 )$ is empty,
\item{C.} $f_1 (D_1 )\cap f_2 (D_2 )$ is non-empty and the 
connecting
tube has its two critical points at the same level. 

See Figure 9.

\centinsert{\pscaption{\psboxscaled{500}{fig6}}{Figure 9}}
\bigskip

 We then define $f:D_1\sharp D_2\to\R$ to be the restriction
 of the height function. Clearly $d f$ blends with $\o_1$ and $\o_2$ 
to
 define a closed $1$-form $\o$ on $M_1\sharp M_2$. Note that in 
each of these 
constructions two new critical points $x,y$ of index 1 are introduced.
In construction A, $f(x)<f(y)$; in construction B,
$f(x)>f(y)$; in construction C,  $f(x)=f(y)$.

If we assume that
 the $\o_i$ are Calabi then these constructions will have two 
important
 properties:
\roster
\item Under construction A, $\o$ is Calabi; under constructions B 
and C $\o$ the consructed form is
 not Calabi. 
\item Under construction A, every leaf of $\o$ intersects at least 
one leaf
 of $\o_1$ and one leaf of $\o_2$ ; under construction B there are 
new
 compact leaves produced; under construction C there is just one 
new compact
singular leaf component produced.

\endroster

We now use these general constructions to produce the four 
examples promised.
 In these
examples we denote the torus by $T$ and $\theta ,\phi$ will denote 
the
 usual angle 
coordinates on $T$.

{\bf Example 1.}    We construct a Calabi form with some compact 
leaves whose cohomology class has rank one less than the first Betti 
number. Let $M_1 =T$ and $\o_1 =d\theta $ (or even 
$pd\theta +
qd\phi$ for $p,q$ relatively prime integers). Let $M_2 ,\o_2$ be
 arbitrary -- for
example $\o_2$ might be a completely irrational Calabi form. 
Choose $D_1\subset M_1$ 
to be small
enough so that many of the leaves of $\o_1$ do not meet $D_1$. If 
we use
 construction A, 
the result is a Calabi form of non-maximal rank with some compact 
leaves.
See Figure 10.

\qquad\qquad\qquad
\centinsert{\pscaption{\psboxscaled{500}{fig7}}{Figure 10}}
\bigskip

{\bf Example 2.}    We construct a Calabi form with all 
non-compact leaves whose cohomology class is split but has rank $>1$. 
Let $M_1 =M_2 =T$ and $\o_1 =d\theta
 ,\o_2 =\lambda d\theta$, where
$\lambda$ is irrational, and we use construction A to produce a 
Calabi form.
 If the disks $D_1$ and $D_2$
are small, then every leaf of $\o_1$ or $\o_2$ will hit
them at most once and the leaves of $\o$ will all be compact. If we 
choose
the disks differently, though, so that they take the form of thin 
ribbons
 winding several times around $T$, then each leaf will hit the disk
 several times. 
The increments $\Delta f_1$ between different intersections of the 
same
 leaf of $\o_1$ with $D_1$ will be integral multiples of $2\pi$.
 The increments $\Delta f_2$ between different intersections of the 
same
 leaf of $\o_2$ with $D_2$ will be integral multiples of 
$2\pi/\lambda$.
 Since $\lambda$ is irrational, the attachments of leaves of $\o_1$ 
with
 leaves of $\o_2$ will produce non-compact leaves of $\o$. See 
Figure 11.  

\qquad\qquad\qquad\qquad
\centinsert{\pscaption{\psboxscaled{500}{fig8}}{Figure 11}}
\bigskip

{\bf Example 3.}    We construct a non-Calabi form with some 
compact leaves whose cohomology class is completely irrational. 
Let $M_i$ and $\o_i$ be arbitrary and use 
construction B.
 This will produce a non-Calabi Morse $1$-form $\o$ with some 
compact
 leaves. In this way choosing $\o_i$'s to be completely irrational 
with corresponding normalization so that $\o$ will be completely
irrational as well, we can realize any cohomology class on any 
surface except for a completely irrational class on a torus.  

{\bf Example 4.}    We construct a non-Calabi non-generic form 
with all non-compact leaves. This is similar to example 3. Let $M_1 
=M_2 =T$ and
$\o_i$ be any two forms of the form $a\,d\theta +b\,d\phi$, where
$a/b$ is irrational. Then use construction C with small disks $D_1$ 
and $D_2$. This will produce a
non-Calabi form with one compact singular leaf component and all 
non-compact leaves. See Figure 12.  

\qquad\qquad\qquad\qquad
\centinsert{\pscaption{\psboxscaled{500}{fig9}}{Figure 12}}
\bigskip

{\bf 13. Final Remarks.} Theorems 1 and 2 still leave us with the 
problem of deciding, in general, whether closed Morse 1-forms are 
contiguous to harmonic forms. 
\footnote{This has recently been settled in the affirmative by Ko Honda.}
The next case to consider might be 
when the form $\o$ has all compact leaves -- the issue then reduces 
to a combinatorial question about the graph $\Gamma_{\o}$. Settling 
the general compact leaf case and using Theorem 2 might then 
enable one to use Proposition 1 to deal with the general case.

The results of Sections 10-11 and the examples of Section 12 leave 
open several interesting questions about the relation between 
the cohomology class of a closed 
Morse 1-form and the presence of compact leaves in its foliation.
\item{(1)} What conditions on a cohomology class will assure that it 
has a
 representative Morse form with no compact leaves? Example (2) 
suggests that
 rank $>1$ might be sufficient.
\item{(2)} What conditions on a cohomology class are satisfied if 
{\it every }
 representative Morse form has all leaves compact? Example (2) 
suggests that
 rank $=1$ is necessary.
\item{(3)} Is there some condition on a cohomology class which 
implies that 
{\it every} representative Morse form has all non-compact leaves? 
It is not
 hard to see that, on the torus, rank $>1$ will do. Example (1) shows 
that 
this is the only case for surfaces. In higher dimensions it is not hard 
to
 see that if a cohomology class $\alpha$ is completely irrational and 
unsplittable, in the sense that it is impossible to separate the 
manifold
 $M$ by a codimension one submanifold $V$ so that $\alpha\not= 0$ 
on each 
component of $M-V$, and if the pairing 
$H^1 (M)\times H^1 (M)\to H^2 (M)$ is non-degenerate then $\alpha$ 
has the
 desired property.

\enddocument